\documentclass[14pt]{article}
\usepackage{amsmath,amssymb}
\renewcommand{\title}[1]{\hrule\begin{center}{\bfseries #1}\par}
\renewcommand{\author}[1]{\medskip{#1}\par\smallskip}
\newcommand{\affiliation}[1]{{\itshape #1}\par}
\newcommand{\email}[1]{e-mail:~\texttt{#1}\par}
\renewcommand{\maketitle}{\end{center}\hrule\par\bigskip}
\renewcommand{\abstract}{\relax}

\begin{document}
%%%%%%%%%%%%%%%%%%%%%%%%%%%%%%%%%%%%%%%%%%%%%%
% Please do not alter the layout!
%%%%%%%%%%%%%%%%%%%%%%%%%%%%%%%%%%%%%%%%%%%%%%
%
%% Title of your talk
\title{\bf Quantum MIMO n-Systems and Conditions for Stability}
 % required
%
%% First author
\author{Seyed M.H. Mansourbeigi }%
\affiliation{Dept.of Electrical Engineering , Polytechnic Univ. Long Island Graduate Center, New York, USA}%
\email{s.mansourbeigi@ieee.org}%
%
% Additional affiliations in the same form...

%% Second author
\author{Vida Milani} % required
\affiliation{Dept.of Math., Faculty of math.Sciences, Shahid Beheshti Univ., Tehran, 1983963113, IRAN } % required
\email{v-milani@cc.sbu.ac.ir} % optional

\maketitle
%%%%%%%%%%%%%%%%%%%%%%%%%%%%%%%%%%%%%%%%%%%%%%

\begin{center}
{\it \bf Abstract}
\end{center}

 … In this paper we present some conditions
for the (strong) stabilizability of an n-D Quantum MIMO system
P(X). It contains two parts. The first part is to introduce the
n-D Quantum MIMO systems where the coefficients vary in the
algebra of Q-meromorphic functions. Then we introduce some
conditions for the stabilizability of these systems. The second
part is to show that this Quantum system has the n-D system as
its quantum limit and the results for the SISO,SIMO,MISO,MIMO are
obtained again as special cases.
\\
\\
\\

\begin{center}
{\bf \large Introduction}
\end{center}

In modern communication technic, multiple-input and
multiple-output, or MIMO, is the use of multiple antennas at both
the receiver and transmitter to improve communication performance.
It is one of several forms of smart antenna (SA), and the state
of the art of SA technology [1], [5].

MIMO technology has attracted attention, since it offers
significant increases in data throughput and link range without
additional bandwidth or transmit power. It achieves this by
higher spectral efficiency (more bits per second Hertz of
bandwodth) and link reliability or diversity (reduced fading).
Because of these properties, MIMO is a current direction of
international wireless research [4].

 The goal in this work is to study the conditions for
stability of n-D systems. A dynamical system can be interpreted
as a vector field on $R^n$ presented by an ordinary differential
equation. In dynamical systems we deal with trajectories. Since
some systems are very complicated, the understanding of
trajectories will be with difficulties and so sometimes in
studying the system, the notion of stability has to be introduced.

In this work we would like to study dynamical systems on the
Quantized spaces. The notion of functional quantization is
introduced in [1].

The paper is organized as follows: First we introduce the notion
of trajectories and vector fields on quantized spaces and obtain
some facts about them which makes us ready to define dynamics of
these systems. Then we see what can be meant by stability of such
systems and study the conditions of stable systems. In the end we
notice that the quantum limit of the q-system is the dynamical
system on the classical space. We end the paper by presenting
some examples of such systems and find their stability conditions
and also their quantum limit.
\\
\\
\\

\begin{center}

{\bf \Large 1. Quantum Trajectories}

\end{center}

Before going further lets recall the definition of functional
quantization and quantized spaces from[1].\\

Let $S_1$ and $S_2$ be sets. Let $x \in S_1$. Assume that $A$ is a
unital $C$-algebra of complex valued functions on $S_1$. Suppose
that $\ss$ is the $C$-vector space of the algebra of complex
valued functions on $S_1 \times S_2$. Let a composition law
$\star$ makes $\ss$ into an associative, unital, not necessary
commutative $A$-algebra. Denote the restriction of an element $f
\in \ss $ to $\{x\} \times S_2$ by $f_x$ and let $$B=\{f_x | f
\in \ss \}$$. This can be considered as the subalgebra of complex
valued functions on $S_2$. Also let define $\delta _x : A
\rightarrow C
$ defined by $\delta_x (g) = g(x)$ be the character of $A$.\\

With these notations $\ss$ is called a $(x, S_1, A)$ {\it
functional quantization} of $B$ and the homomorphism $$\phi : \ss
\rightarrow B$$ defined by $\phi (f) = f_x$ with the property
$\phi (gf) = \delta_x(g) \phi (f)$ for all $f \in \ss $ and $g
\in A $, is
called the {\it quantization map}.\\

By a quantum space we mean the functional quantization of the
coordinate algebra of $R^n$ as follows:\\

For the coordinates $x$ and $y$ in $R^2$, let $\ast$ be defined as
$y \ast x = q x \ast y$ with $ q \in D-\{0\}$.\\

Now let $D=\{q \in C| |q| \leq 1 \}$ be the unit disc in $C$ and
$A_1(q)$ be the $C$-algebra of all absolutely convergent power
series $\sum_{i=0}^\infty a_i q^i$ in $D$ with coefficients in
$C$. Also let $A_0(q)$ be the $C$-algebra of all absolutely
convergent power series $\sum_{i>- \infty} c_i q^i$ in $D - \{0\}$
with coefficients in $C$. The $( 1, D-\{0\}, A_0(q))$ functional
quantization of $M$: the $C$-algebra of all absolutely convergent
power series $\sum_{i>>- \infty} a_{ij}t_1^i t_2^j$ on
$R-\{0\}\times R-\{0\}$ with coefficients in $C$, is called the
{\it quantum 2-space} and we denote it by $Q^2$. This is a unital
non commutative associative $A_0(q)$-algebra. Elements of this
algebra has a representation of the form $\sum_{i>>- \infty}
a_{ij}x^i
y^j$, where $a_{ij}$s are in $A_0(q)$.\\

{\bf Definition 1.1}. Let $\psi : A_0(q) \rightarrow R $ be the
$C$-algebra homomorphism defined by $\psi (q) = 1$. Let $I
\subseteq R$ be an open interval. Any $\psi$- homomorphism
$\alpha : Q^2 \rightarrow C^{\infty}(I)$ is called a {\it quantum
trajectory} on the quantum 2-space $Q^2$. We call it in short a
Q-trajectory. In other words a Q-trajectory is a $C$-linear map
$\alpha $ satisfying $\alpha (q^i x^j y^k) = \psi (q)^i \alpha
(x)^j \alpha (y)^k $.\\

{\bf Remark }. From the above property we see that $\alpha $ is
completely determined by its values $x$ and $y$, i.e by $\alpha
(x)$
and $\alpha (y)$.\\

{\bf Definition 2.1}. For the point $P=(p_1,p_2) \in R^2$ we say
that $\alpha $ passes through $P$ at $t=0$ if for each $f \in
Q^2$, $ \alpha (f)(0) = f(P)$. In particular if $\alpha $ passes
through $P$ at $t=0$, then $\alpha (x)(0) = p_1$ and $\alpha
(y)(0) = p_2$. And so we can use the familiar notation $\alpha
(0) = (\alpha (x)(0) ,
\alpha (y) (0)) = (p_1, p_2) = P$.\\

{\bf Definition 3.1}. If we write $\alpha $ as $\alpha (t) =
(\alpha (x)(t) , \alpha (y)(t))$, then the {\it velocity vector}
for $\alpha $ at $t$ is defined by $$ \dot{\alpha}(t) = (
\frac{d}{dt} \alpha (x)(t) , \frac{d}{dt} \alpha (y)(t))$$\\

{\bf Definition 4.1}. Let $f \in Q^2$ and let $\alpha : Q^2
\rightarrow C^\infty (I) $ be a Q-trajectory. The rate of change
of $f$ in $\alpha $ direction $\dot{\alpha}f$ is defined by
$$\dot{\alpha }f (t)= \frac{d}{dt} \alpha (f)(t)$$\\

{\bf Lemma 5.1}. $\dot {\alpha}(t)$ can be considered as a
derivation of $Q^2$, i.e. as a $C$-linear map $\dot {\alpha} : Q^2
\rightarrow C^\infty (I)$ with the property that for each $f$ and
$g$ in $Q^2$: $$\dot{\alpha }(fg) (t)= \dot{\alpha }(f)(t).
\dot{\alpha}(g)(t) + \alpha (f)(t) \dot{\alpha }(g)(t)$$\\

{\bf Definition 6.1}. Any $A_0(q)$-homomorphism $X : Q^2
\rightarrow Q^2$ is called a {\it vector field} on $Q^2$. A
Q-trajectory $\alpha$ is the integral curve for $X$ if $$\alpha o
X (t) = \frac{d}{dt}\alpha (t)$$\\

{\bf Lemma 7.1}. For vector fields $X$ and $Y$ their composition
is a vector field and so if we define their bracket by $$[X,Y] =
XoY - YoX$$ then the set of vector fields with this bracket is a
Lie algebra. Furthermore the Liebniz rule satisfies, i.e. for
vector fields $X$ , $Y$ and $Z$ $$[X , YZ] = [X , Y] Z + Y [X ,
Z]$$ Where we define the product of two vector fields by $$XY (f)
= X(f). Y(f)$$ for all $f \in Q^2$.\\

{\bf Proposition 8.1}. For each vector field $X: Q^2 \rightarrow
Q^2$ and each point $P=(p_1,p_2) \in R^2$, there exists a
Q-trajectory $\alpha : Q^2 \rightarrow C^\infty (I)$ passing
through $P$ and satisfies $$\alpha o X (0) = \dot{\alpha} (0)$$\\

{\bf Proof}: Let $X(x) = f $ and $X(y) = g$. Set $$\alpha (t) = (
\alpha (x) (t) , \alpha (y)(t) ) = (p_1 + t \psi (f(p_1,p_2)) ,
P_2 + t \psi (g(p_1,p_2)))$$ Then we have $\alpha (0) = (\alpha
(x)(0), \alpha (y)(0)) = (p_1 , p_2)$ and
$$ \dot{\alpha}(x)(0) = \psi (f(p_1,p_2)) = \alpha (f)(0) = (\alpha o X) (x)(0)$$
$$ \dot{\alpha}(y)(0) = \psi (g(p_1,p_2)) = \alpha (g)(0) = (\alpha o X)
(y)(0)$$\\

{\bf Remark}. The converse of the above proposition is also true.
That is for any Q-trajectory $\alpha : Q^2 \rightarrow C^\infty
(I)$ if $\alpha (0) = (\alpha (x)(0) , \alpha (y)(0)) = (p_1 ,
p_2) = P$, then there exists a vector field $X : Q^2 \rightarrow
Q^2$ satisfying $$\alpha o X (x)(0) = \dot{\alpha} (x)(0)$$
$$\alpha o X (y)(0) = \dot{\alpha} (y)(0)$$\\

{\bf Proof}: It is sufficient to define $X(x)= f$ and $X(y)= g$
where $$\psi (f(p_1,p_2)) = \dot{\alpha} (x)(0)$$
$$\psi (g(p_1,p_2)) = \dot{\alpha} (y)(0)$$\\

\begin{center}

{\bf \Large 2. Quantum Systems}\\

\end{center}

An autonomous (vector field) equation of the form $$P(X): \alpha
o X (t) = \frac{d}{dt} \alpha (t)$$ where $\alpha$ is a
Q-trajectory and $X : Q^2 \rightarrow Q^2 $ is a vector field, is
called an {\it autonomous quantum (MIMO) system}. A solution for
this system is a Q-trajectory $\alpha (t)$ such that $\alpha o X
(t) = 0$.

A {\it nonautonomous quantum (MIMO) system} is defined the same
way except that any Q-trajectory $\alpha$ satisfying the equation
is
a solution for the system.\\

{\bf Example 1.2}. Let $X : Q^2 \rightarrow Q^2$ be the vector
field defined by $X(x) = y$ and $X(y) = x$. Then an autonomous
(time independent) solution for this system is the Q-trajectory
$\alpha$ defined by $\alpha (t) = (0,0)$ and a nonautonomous
(time dependent) solution is the Q-trajectory
$$\alpha (t) =
(\alpha (x)(t) , \alpha (y)(t)) = (c e^{-t} , -c e^{-t})$$\\

{\bf Definition 2.2}. The solution $\alpha$ is called {\it stable}
if for every $\epsilon > 0$ there exists a $\delta(\epsilon) > 0$
such that if $\beta$ is another solution for the system and if
$$|\alpha(t_0)(x) - \beta(t_0)(x)| < \delta ,
|\alpha(t_0)(y) - \beta(t_0)(y)| < \delta$$ then $$|\alpha(t)(x)
- \beta(t)(x)| < \epsilon , |\alpha(t)(y) - \beta(t)(y)| <
\epsilon$$ for each $t > t_0$ and $t_0 \in I$.\\

{\bf Example 3.2}. Both the nonautonomous and autonomous solutions
for example 1.2 are stable.\\

{\bf Note}. We can generalize all the above results for the
n-dimensional case.\\

\begin{center}

{\bf \Large 3. Quantum Limit}\\

\end{center}

In the quantum limit $$q \rightarrow 0$$ the equations of motion
and the evolution of the quantum system will have the classic
meaning in terms of the classical Hamiltonian.

\end{document}